\documentclass[11pt,twoside]{amsart}
\usepackage{latexsym,amssymb,amsmath}

\numberwithin{equation}{section}
\newtheorem{theorem}{Theorem}[section]
\newtheorem{lemma}[theorem]{Lemma}
\newtheorem{proposition}[theorem]{Proposition}
\newtheorem{corollary}[theorem]{Corollary}

\theoremstyle{definition}
\newtheorem{definition}[theorem]{Definition} 
\newtheorem{remark}[theorem]{Remark}
\newtheorem{example}[theorem]{Example}
\newtheorem{problem}[theorem]{Problem}

\begin{document}

\title[Duality, a-invariants and canonical modules]
{Duality, a-invariants and canonical modules         
of rings arising from linear optimization problems
}
 
\author{Joseph P. Brennan}
\address{Department of Mathematics \\
University of Central Florida\\
Orlando, FL 32816-1364, USA.}
\email{jpbrenna@mail.ucf.edu}

\author{Luis A. Dupont}
\address{
Departamento de
Matem\'aticas\\
Centro de Investigaci\'on y de Estudios
Avanzados del
IPN\\
Apartado Postal
14--740 \\
07000 Mexico City, D.F.
}
\email{ldupont@math.cinvestav.mx}

\author{Rafael H. Villarreal}
\address{
Departamento de
Matem\'aticas\\
Centro de Investigaci\'on y de Estudios
Avanzados del
IPN\\
Apartado Postal
14--740 \\
07000 Mexico City, D.F.
}
\email{vila@math.cinvestav.mx}

  
\begin{abstract} The aim of this paper is to study integer
rounding properties of various systems of linear inequalities to
gain insight about  
the algebraic properties of Rees algebras of monomial ideals and
monomial subrings. We study the normality and Gorenstein
property---as well as 
the canonical module and the $a$-invariant---of Rees
algebras and subrings arising from systems with
the integer rounding property. We relate the algebraic properties of
Rees algebras and monomial subrings with integer rounding properties
and present a duality theorem. 
\end{abstract}
 
\maketitle

\section{Introduction}

Let $R=K[x_1,\ldots,x_n]$ be a polynomial ring over a field $K$ and
let $v_1,\ldots,v_q$ be the column vectors of a matrix $A=(a_{ij})$
whose entries are 
non-negative integers. We shall always assume that the rows and
columns of $A$ are different from zero. As usual we
use the notation $x^a:=x_1^{a_1} \cdots x_n^{a_n}$, 
where $a=(a_1,\ldots,a_n)\in \mathbb{N}^n$. 

The {\it monomial algebras\/} considered here are: 
(a) the {\it Rees algebra\/}
$$
R[It]:=R\oplus It\oplus\cdots\oplus I^{i}t^i\oplus\cdots
\subset R[t],
$$
where $I=(x^{v_1},\ldots,x^{v_q})\subset R$ and $t$ is a new
variable, (b) the {\it extended Rees algebra\/} 
$$
R[It,t^{-1}]:=R[It][t^{-1}]\subset R[t,t^{-1}],
$$
(c) the 
{\it monomial subring\/} 
$$
K[F]=K[x^{v_1},\ldots,x^{v_q}]\subset R
$$
spanned by $F=\{x^{v_1},\ldots,x^{v_q}\}$, (d) the 
{\it homogeneous monomial subring\/} 
$$
K[Ft]=K[x^{v_1}t,\ldots,x^{v_q}t]\subset R[t]
$$
spanned by $Ft$, (e) the 
{\it homogeneous monomial subring\/} 
$$
K[Ft\cup\{t\}]=K[x^{v_1}t,\ldots,x^{v_q}t,t]\subset R[t]
$$
spanned by $Ft\cup\{t\}$, (f) the {\it homogeneous
monomial subring\/}
$$
S=K[x^{w_1}t,\ldots,x^{w_r}t]\subset R[t],
$$        
where $w_1,\ldots,w_r$ is the set of all vectors
$\alpha\in\mathbb{N}^n$ such 
that $0\leq \alpha\leq v_i$ for some $i$, and (g) the {\it Ehrhart
ring\/}  
$$
A(P)=K[\{x^at^i\vert\, a\in \mathbb{Z}^n \cap iP; i\in
\mathbb{N}\}]\subset R[t] 
$$
of a lattice polytope $P$. 

The aim of this work is to study max-flow min-cut
properties of clutters and integer rounding properties of
various systems of linear inequalities---and their underlying
polyhedra---to gain insight about the algebraic properties of these 
algebras and viceversa. Systems with integer rounding properties and
clutters with the max-flow min-cut property come from linear
optimization problems \cite{Schr,Schr2}. The precise 
definitions will be given in Section~\ref{i-r-p}.

Before stating our main results, we recall a few basic facts about the
normality of monomial subrings. 
According to \cite{monalg} the 
{\it integral closure\/} of $K[F]$ 
in its field of fractions can be expressed as
\begin{equation}\label{intclos-desc-08}
\overline{K[F]}=K[\{x^a\vert\, a\in
\mathbb{Z}{\mathcal
A}\cap \mathbb{R}_+{\mathcal A}\}],
\end{equation}
where $\mathcal{A}=\{v_1,\ldots,v_q\}$, $\mathbb{Z}{\mathcal A}$ is the
subgroup of $\mathbb{Z}^n$
spanned by ${\mathcal A}$, and $\mathbb{R}_+\mathcal{A}$ is the cone
generated 
by $\mathcal{A}$.  The subring $K[F]$ equals
$K[\mathbb{N}\mathcal{A}]$, the semigroup ring of 
$\mathbb{N}\mathcal{A}$. Recall that $K[F]$ is called {\it integrally
closed\/} or {\it normal\/} if $K[F]=\overline{K[F]}$. 
Thus $K[F]$ is normal if and
only if  
$$
\mathbb{N}\mathcal{A}=\mathbb{Z}\mathcal{A}\cap\mathbb{R}_+{\mathcal A},
$$
where $\mathbb{N}\mathcal{A}$ is the subsemigroup of $\mathbb{N}^n$
generated by $\mathcal{A}$. The description of the integral
closure given in Eq.~(\ref{intclos-desc-08}) can of course be applied
to any of the monomial algebras considered here. In particular if 
$\mathcal{A}'$ is the set
$$
\mathcal{A}'=\{e_1,\ldots,e_n,(v_1,1),\ldots,(v_q,1)\},
$$
where $e_{i}$ is the $i${\it th} unit vector, then
$\mathbb{Z}\mathcal{A}'=\mathbb{Z}^{n+1}$ and $R[It]$ is normal if and
only if 
$\mathbb{N}\mathcal{A}'=\mathbb{Z}^{n+1}\cap\mathbb{R}_+{\mathcal
A}'$. A dual characterization of the normality of $R[It]$ will be given
in Proposition~\ref{jun25-08}. 

Recall that the Ehrhart ring $A(P)$ is always normal \cite{BHer}.  A set 
$\mathcal{A}\subset\mathbb{Z}^n$ is called a {\it Hilbert basis\/} if
$\mathbb{N}\mathcal{A}=\mathbb{R}_+\mathcal{A}\cap\mathbb{Z}^n$. Note
that if $\mathcal{A}$ is a Hilbert basis, then the ring $K[F]$ is normal. 

The contents of this paper are as
follows. First we use the theory of blocking and antiblocking polyhedra
\cite{baum-trotter,fulkerson,fulkerson-jctb,Schr} to describe when
the systems  
$$x\geq 0;\,
xA\leq\mathbf{1},\ \ \  x\geq 0;xA\geq\mathbf{1},\ \ \
xA\leq\mathbf{1},
$$
have the {\it integer rounding 
property\/} (see Definitions~\ref{irpx<=1}, \ref{irpx>=1}, \ref{irp})
in terms of 
the normality of the monomial algebras considered here. As usual, 
we denote the vector $(1,\ldots,1)$ by $\mathbf{1}$. 
If $a=(a_1,\ldots,a_n)$ and $b=(b_1,\ldots,b_n)$ are vectors, we 
write $a\leq b$ if $a_i\leq b_i$ for all $i$.

One of the main results of Section~\ref{i-r-p} is:

\medskip 

\noindent {\bf Theorem~\ref{round-up-char}}{\it\ The system $x\geq
0;\, xA\leq \mathbf{1}$ has the integer 
rounding property if and only if the subring
$S=K[x^{w_1}t,\ldots,x^{w_r}t]$ is normal.} 

\medskip

This result was shown in \cite{roundp} when $A$ is the incidence
matrix of a clutter, i.e., when the entries of $A$ are in $\{0,1\}$. 
Recall that a {\it clutter\/} $\mathcal{C}$ with
finite vertex set $X=\{x_1,\ldots,x_n\}$ 
is a family of subsets of $X$, called edges, none 
of which is included in another. The {\it incidence matrix\/} of
a clutter $\mathcal{C}$ is the vertex-edge matrix whose columns are
the  
characteristic vectors of the edges of $\mathcal{C}$. 
The {\it edge ideal\/} of a clutter $\mathcal{C}$,
denoted by $I(\mathcal{C})$, is the ideal of $R$ generated by 
all monomials $x_e=\prod_{x_i\in e}x_i$ such that $e$ is an edge of
$\mathcal{C}$. The {\it Alexander dual\/} of $I(\mathcal{C})$ is the
ideal of $R$ given 
by $I(\mathcal{C})^\vee=\cap_{e\in
E}(e)$, where $E=E(\mathcal{C})$ is the edge set of
$\mathcal{C}$.

The integer rounding property of some systems has already been
expressed in terms of the 
normality of monomial algebras \cite{poset,roundp}. In \cite{poset} 
it is shown that  
the system $x\geq 0; xA\geq \mathbf{1}$ has the integer rounding
property if 
and only if $R[It]$ is normal (this was also observed by N. V. Trung 
if $A$ 
is the incidence matrix of a clutter). Here we complement this 
fact by presenting a duality between the integer rounding property of
the   systems $x\geq 0; xA\geq \mathbf{1}$ and $x\geq 0; xA^*\leq
\mathbf{1}$ valid for matrices with entries in $\{0,1\}$, where
$a_{ij}^*=1-a_{ij}$ is the $ij$-entry of $A^*$. This duality
is extended to a duality between monomial subrings. 

Altogether 
another main result of Section~\ref{i-r-p} is:

\medskip

\noindent {\bf Theorem~\ref{duality-subrings-irp}}{\it\ Let $A$ be
the incidence 
matrix of a clutter. If $v_i^*=\mathbf{1}-v_i$ and $A^*$ is the
matrix with column vectors 
$v_1^*,\ldots,v_q^*$, then the following are equivalent:
\begin{itemize}
\item[(a)] $R[I^*t]$ is normal, where $I^*=(x^{v_1^*},\ldots,x^{v_q^*})$. 
\item[(b)] $S=K[x^{w_1}t,\ldots,x^{w_r}t]$ is normal.
\item[(c)] $\{-e_1,\ldots,-e_n,(v_1,1),\ldots,(v_q,1)\}$ is a
Hilbert basis. 
\item[(d)] $x\geq 0;xA^*\geq\mathbf{1}$ has the integer rounding
property.
\item[(e)] $x\geq 0;xA\leq\mathbf{1}$ has the integer rounding
property.
\end{itemize}}

Then we present some interesting consequences of this duality. First
of all we recover one of the main results of \cite{perfect} showing
that if 
$$P=\{x\vert x\geq 0;xA\leq\mathbf{1}\}$$
is an 
integral polytope, i.e., $P$ has only integer vertices, 
and $A$ is a $\{0,1\}$-matrix, then the Rees algebra
$R[I^*t]$ is normal (see Corollary~\ref{jun26-08}). This result is
related to 
perfect graphs. Indeed if $P$ is
integral, then $v_1,\ldots,v_q$ correspond to the maximal cliques
(maximal  
complete subgraphs) of a perfect graph $H$ \cite{chvatal,lovasz}, and 
$v_1^*,\ldots,v_q^*$ correspond to the minimal vertex covers of the
complement of $H$. Second
we show that if $A$ is the incidence matrix of the collection of basis
of a matroid, then all systems 
$$
x\geq 0;xA\geq\mathbf{1},\ \ \ x\geq 0;xA^*\geq\mathbf{1},\ \ \ 
x\geq 0;xA\leq\mathbf{1},\ \ \ x\geq 0;xA^*\leq\mathbf{1}
$$
have the integer rounding property (see Corollary~\ref{jun29-08}). 
Third we show that if $A$ is the incidence matrix of a graph, then 
$R[It]$ is normal if and only if $R[I^*t]$ is normal (see
Corollary~\ref{forgraphs-In-iff-dualIn}). 
We give an example to show that this result does not extends to
arbitrary uniform clutters (see
Example~\ref{forclutters-In-iff-dualIn-false}). If $A$ is the
incidence matrix of a graph $G$, we characterize when $I^*$ is the
Alexander dual of the edge ideal of the complement of $G$ (see 
Proposition~\ref{el-vegetariano}). If $G$ is a triangle-free graph,
we show a 
duality between the normality of $I=I(G)$ and that of the Alexander
dual  of the edge ideal of the complement of $G$ (see
Corollary~\ref{july4-08}).  
We show an example of an edge 
ideal of a graph whose Alexander dual is not normal (see
Example~\ref{contraec}). In \cite{perfect} it is shown
that this is never the case if the graph is perfect, i.e., the
Alexander dual of the edge ideal of a perfect graph is always normal.  
Finally we recover one of the main results of \cite{clutters} showing
that if $A$ is the incidence matrix of a clutter $\mathcal{C}$, then 
$\mathcal{C}$ satisfies the max-flow min-cut
property if and only if the set covering polyhedron 
$$Q(A)=\{x\vert\, x\geq
0;xA\geq\mathbf{1}\}$$ 
is integral and $R[It]$ is normal (see
Corollary~\ref{jun29-08-1}). 

The last main result of
Section~\ref{i-r-p} is:  

\medskip

\noindent {\bf Theorem~\ref{jun5-08}}{\it\ If the system
$xA\leq\mathbf{1}$ has the integer rounding property, then $K[F]$ 
is normal and $\mathbb{Z}^n/\mathbb{Z}\mathcal{A}$ is a torsion-free
group. 
The converse holds if $|v_i|=d$ for all $i$. Here
$|v_i|=\langle v_i,\mathbf{1}\rangle$}.

\medskip

As a consequence of this result we prove: (i) If $A$ is the incidence
matrix of a  
connected graph $G$, then the system $xA\leq\mathbf{1}$ has the
integer rounding property if and only if $G$ is a bipartite graph 
(see Corollary~\ref{jun8-08-1}), and (ii)  Let $A$ be the incidence
matrix of a clutter 
$\mathcal{C}$. If $\mathcal{C}$ is uniform, i.e., all its edges have the
same size,  and $\mathcal{C}$ has the max-flow min-cut
property (see Definition~\ref{mfmc-def}), then the system
$xA\leq\mathbf{1}$ has the 
integer rounding property (see Corollary~\ref{jun18-08}). 

If $A$ is the incidence matrix of a
bipartite graph, a remarkable result of \cite{roundp} shows that 
the system $x\geq 0; xA\leq\mathbf{1}$ has the integer rounding
property if and only if the extended Rees algebra $R[It,t^{-1}]$ is
normal. 

Before stating the main results
of Sections~\ref{can-mod-device} and \ref{section-on-canmod}, 
we need to introduce 
the canonical module and the $a$-invariant (see
Section~\ref{can-mod-device} for additional details). Below we
briefly explain the 
important role that these two objects play in the general theory. The
subring $S$ is a standard $K$-algebra because 
$\langle(w_i,1),e_{n+1}\rangle=1$ for all $i$. Here $\langle\,,\rangle$
is the standard inner product and $e_{i}$ is the $i${\it th} unit
vector. If $S$ is normal, then according to a formula of Danilov and 
Stanley \cite{Dan} the canonical module of $S$ is 
the ideal of $S$ given by  
\begin{equation}
\omega_{S}=(\{x^at^b\vert\,
(a,b)\in \mathbb{N}{\mathcal B}\cap ({\mathbb R}_+{\mathcal B})^{\rm
o}\}), 
\end{equation}
where  ${\mathcal B}=\{(w_1,1),\ldots,(w_r,1)\}$ and 
$({\mathbb R}_+{\mathcal B})^{\rm o}$ is the relative interior of 
${\mathbb R}_+{\mathcal B}$. This expression for the canonical module
of $S$ is 
central for our purposes. Recall that the $a$-{\it invariant\/} of
$S$, denoted by $a(S)$, is  
the degree as a rational 
function of the Hilbert series 
of $S$ \cite[p.~99]{monalg}. Thus we may compute $a$-invariants using
the program {\it Normaliz} \cite{B}. Let $H_S$ and
$\varphi_S$ be the Hilbert function and the Hilbert polynomial of $S$
respectively. The {\it index of regularity\/} of $S$, denoted by ${\rm
reg}(S)$, is the least positive integer such that $H_S(i)=\varphi_S(i)$
for $i\geq {\rm reg}(S)$. The $a$-invariant plays a fundamental role
in algebra and geometry because one has: 
$$
{\rm reg}(S)=\left\{\begin{array}{ll}
0&\mbox{if }a(S)<0,\\
a(S)+1&\mbox{otherwise},
\end{array}
\right.
$$
see \cite[Corollary~4.1.12]{monalg}. If $S$ is 
normal, then $S$ is Cohen-Macaulay \cite{Ho1} and its $a$-invariant is
given by
\begin{equation}\label{princeton-fall-07-1}
a(S)=-{\rm min}\{\, i\, \vert\, (\omega_S)_i\neq 0\},
\end{equation}
see \cite[p.~141]{BHer} and \cite[Proposition~4.2.3]{monalg}. 

In Section~\ref{can-mod-device} we give a general technique to
compute the 
canonical module and
the $a$-invariant of a wide class of monomial subrings (see
Theorem~\ref{canmod-tdi}). 

Then in Section~\ref{section-on-canmod} we
study the canonical module and 
the $a$-invariant of monomial subrings arising from integer rounding
properties. We give necessary and sufficient conditions for $S$ to be
Gorenstein and give a formula for the $a$-invariant of $S$ in terms 
of the vertices of the polytope $P=\{x\vert\, x\geq 0; xA\leq
\mathbf{1}\}$. For use below let ${\rm vert}(P)$ be the set
of vertices of $P$  and 
let $\ell_1,\ldots,\ell_p$ be the set of all maximal elements of ${\rm
vert}(P)$ (maximal with respect to $\leq$). 
For each $1\leq i\leq p$ there is a
unique positive  
integer $d_i$ such that the non-zero entries of $(-d_i\ell_i,d_i)$ are
relatively prime. 

The main results of Section~\ref{section-on-canmod}
are as follows. 

\medskip

\noindent {\bf Theorem~\ref{can-mod-intr-norm}}{\it\   
If the system $x\geq 0; xA\leq \mathbf{1}$ has the 
integer rounding property, then the canonical module of
$S=K[x^{w_1}t,\ldots,x^{w_r}t]$ is given by   
\begin{equation}\label{march5-08-2}
\omega_S=\left(\left\{\left.x^at^b\right\vert\,
(a,b)\left(\begin{array}{rrrlrr} 
-d_1\ell_1&\cdots&-d_p\ell_p&e_1&\cdots& e_n\\
d_1        &\cdots& d_p     &\ 0  &  \cdots&0
\end{array}
\right)\geq\mathbf{1}\right\}\right),
\end{equation}
and the $a$-invariant of $S$ is equal to 
$-\max_i\{\lceil 1/d_i+|\ell_i|\rceil\}$. Here
$|\ell_i|=\langle\ell_i,\mathbf{1}\rangle$.} 

\medskip

This result complements a result of \cite{roundp}
valid only for incidence matrices of clutters. If
$S$ is normal, the last Betti number in the homogeneous free
resolution of the toric ideal $P_S$ of $S$ is equal to
$\nu(\omega_S)$, the minimum number of generators of $\omega_S$. This
number is called the {\it type\/} of $P_S$. 
Thus by describing the canonical module of $S$ we are in fact
providing a device to compute the type of $P_S$. According to
\cite{shiftcon} the number of integral vertices of the polyhedron that
defines $\omega_S$ (see Eq.~(\ref{march5-08-2})) is a lower bound for
$\nu(\omega_S)$.

Using the description above for $\omega_S$ we then prove:

\medskip

\noindent {\bf Theorem~\ref{may27-08-1}}{\it\ Assume that the system
$x\geq 0$; $xA\leq\mathbf{1}$ has  
the integer rounding property. If $S$ is
Gorenstein and $c_0=\max\{|\ell_i|\colon\, 1\leq i\leq p\}$ is an
integer, then $|\ell_k|=c_0$ for each $1\leq k\leq p$ such that
$\ell_k$ has integer entries.}

\medskip

\noindent {\bf Theorem~\ref{may27-08-2}}{\it\  Assume that the system
$x\geq 0; 
xA\leq \mathbf{1}$ 
has the 
integer rounding property. If $-a(S)=1/d_i+|\ell_i|$ for
$i=1,\ldots,p$, then $S$ is 
Gorenstein.}

\medskip

As a consequence of Theorems~\ref{may27-08-1} and \ref{may27-08-2} we
obtain that if $P$ is an integral polytope, i.e., it has only integral
vertices, then $S$ is Gorenstein if
and only if $a(S)=-(|\ell_i|+1)$ for $i=1,\ldots,p$ (see 
Corollary~\ref{may27-08-3}). 

We also examine the Gorenstein and complete intersection 
properties of subrings arising from systems with the integer rounding
property of incidence matrices of graphs. Let $G$ be a connected 
graph with $n$ vertices and $q$ edges and let $A$ be its incidence
matrix. 
Based on a computer
analysis, using the  
program {\it Normaliz\/} \cite{B}, we conjecture a possible description of 
all Gorenstein subrings $S$ in terms of the vertices of $P$ 
(see Problem~\ref{gorenstein-conjecture}).  
If the system $xA\leq\mathbf{1}$ has the integer rounding
property, then we show that $K[Ft\cup\{t\}]$ is 
a complete intersection if and only if $G$ is bipartite and 
the number of primitive cycles of $G$ is equal to $q-n+1$ 
(see Proposition~\ref{jun19-08}).

Let $G$ be a bipartite graph and let $A$ be its
incidence matrix. A constructive description of all bipartite graphs
such that $K[G]=K[x^{v_1},\ldots,x^{v_q}]$ is a complete intersection
is given in \cite{electronic-n}. The Gorenstein property of $K[G]$ has
been studied in \cite{accota-gv,stable-subring}. Thus by
Lemma~\ref{jun8-08} and \cite[Proposition~3.1.19]{BHer} the
Gorenstein property 
and the complete intersection property of
$K[x^{v_1}t,\ldots,x^{v_q}t,t]$ are well understood in this
particular case. The $a$-invariant of $K[G]$ has a combinatorial
expression in terms of 
directed cuts and can be computed using linear programming
\cite{shiftcon}. Some other expressions for $a(K[G])$ can be found in
\cite{corso-nagel-trans,accota-gv,join}.

\section{Integer rounding properties}\label{i-r-p}

We continue to use the notation and definitions used in the
introduction. In 
this section we introduce and study integer rounding properties,
describe some  
of their properties, present a duality theorem and show several
applications. 

Let $P$ be a rational
polyhedron in $\mathbb{R}^n$. Recall that the {\it antiblocking
polyhedron\/} of $P$ is defined as:
$$
T(P):=\{z\vert\, z\geq 0; \langle z,x\rangle\leq 1\mbox{ for all }x\in
P\}. 
$$

\begin{lemma}\label{napkin-lemma}
Let $A$ be a matrix of order $n\times q$ with entries in $\mathbb{N}$, 
let $v_1,\ldots,v_q$ be the column vectors of $A$ and let
$\{w_1,\ldots,w_r\}$ 
be the set of all $\alpha$ in $\mathbb{N}^n$ such that $\alpha\leq v_i$
for some $i$. If $P=\{x\vert\, x\geq 0;\, xA\leq\mathbf{1}\}$, then 
$$
T(P)={\rm conv}(w_1,\ldots,w_r).
$$
\end{lemma}

\begin{proof} First we show the following equality
which is interesting in its own right:
\begin{equation}\label{napkin-equation}
{\rm conv}(w_1,\ldots,w_r)=\mathbb{R}_+^n\cap({\rm
conv}(w_1,\ldots,w_r)+\mathbb{R}_+\{-e_1,\ldots,-e_n\}).
\end{equation}
Clearly the left hand side is contained in the right hand side. 
Conversely let $z$ be a vector in the right hand side. Then $z\geq 0$
and we can
write
\begin{equation}\label{march1-08}
z=\lambda_1w_1+\cdots+\lambda_rw_r-\delta_1e_1-\cdots-\delta_ne_n,\
\ (\lambda_i\geq 0;\, \textstyle\sum_i\lambda_i=1;\, \delta_i\geq 0).
\end{equation}
Consider the vector $z'=\lambda_1w_1+\cdots+\lambda_rw_r-\delta_1e_1$.
We set $T'={\rm conv}(w_1,\ldots,w_r)$ and 
$w_i=(w_{i1},\ldots,w_{in})$.
We claim that $z'$ is in $T'$. We may assume that $\delta_1>0$,
$\lambda_i>0$ for all $i$, and that the first entry $w_{i1}$ of $w_i$
is positive for $1\leq i\leq s$ and is equal to zero for $i>s$. From
Eq.~(\ref{march1-08}) we get $\lambda_1w_{11}+\cdots+\lambda_sw_{s1}\geq
\delta_1$. 

Case (I): $\lambda_1w_{11}\geq \delta_1$. Then we can write
$$
z'=\frac{\delta_1}{w_{11}}\left(w_1-w_{11}e_1\right)+
\left(\lambda_1-\frac{\delta_1}{w_{11}}\right)w_1+
\lambda_2w_2+\cdots+\lambda_rw_r.
$$
Notice that $w_1-w_{11}e_1$ is again in $\{w_1,\ldots,w_r\}$. Thus
$z'$ is a convex combination of $w_1,\ldots,w_r$, i.e.,  
$z'\in T'$.

Case (II): $\lambda_1w_{11}<\delta_1$. Let $m$ be the largest integer less
than or equal to $s$ such that
$\lambda_1w_{11}+\cdots+\lambda_{m-1}w_{(m-1)1}<\delta_1\leq
\lambda_1w_{11}+\cdots+\lambda_{m}w_{m1}$. Then 
\begin{eqnarray*}
z'&=& \sum_{i=1}^{m-1}\lambda_i(w_i-w_{i1}e_1)+
\left[\frac{\delta_1}{w_{m1}}-\left(\sum_{i=1}^{m-1}
\frac{\lambda_iw_{i1}}{w_{m1}}\right)\right](w_m-w_{m1}e_1)+
\\
& & \left[\lambda_m-\frac{\delta_1}{w_{m1}}+\left(\sum_{i=1}^{m-1}
\frac{\lambda_iw_{i1}}{w_{m1}}
\right)\right]w_m+\sum_{i=m+1}^r\lambda_{i}w_{i}.
\end{eqnarray*}
Notice that $w_i-w_{i1}e_1$ is again in $\{w_1,\ldots,w_r\}$ 
for $i=1,\ldots,m$. Thus $z'$ is a convex combination of
$w_1,\ldots,w_r$, i.e., $z'\in T'$. This completes the proof of the
claim. Note that we  
can apply the argument above to any entry of $z$ or $z'$ thus we
obtain that $z'-\delta_2e_2\in T'$. Thus by induction 
we obtain that $z\in T'$, as required. This completes the proof of
Eq.~(\ref{napkin-equation}).

Clearly one has the equality $P=\{z\vert\, z\geq 0;\langle
z,w_i\rangle\leq 1\, \forall i\}$ because for each $w_i$ there is
$v_j$ such that $w_i\leq v_j$. Hence by the finite basis theorem
\cite{Schr} we can write 
\begin{equation}\label{march6-08-1}
P=\{z\vert\, z\geq 0;\langle z,w_i\rangle\leq 1\, \forall i\}={\rm
conv}(\ell_0,\ell_1,\ldots,\ell_m)
\end{equation}
for some $\ell_1,\ldots,\ell_m$ in
$\mathbb{Q}_+^n$ and $\ell_0=0$. From Eq.~(\ref{march6-08-1}) we
readily get the equality
\begin{equation}\label{march6-08-3}
\{z\vert\, z\geq 0;\langle z,\ell_i\rangle\leq 1\, \forall i\}=T(P).
\end{equation}
Using Eq.~(\ref{march6-08-1}) and noticing that 
$\langle\ell_i,w_j\rangle\leq 1$ for all $i,j$, we get
$$
\mathbb{R}_+^n\cap({\rm
conv}(\ell_0,\ldots,\ell_m)+\mathbb{R}_+\{-e_1,\ldots,-e_n\})=
\{z\vert\, z\geq 0;\langle z,w_i\rangle\leq 1\, \forall i\}.
$$
Hence using this equality and \cite[Theorem~9.4]{Schr} we obtain
\begin{equation}\label{march6-08-2}
\mathbb{R}_+^n\cap({\rm
conv}(w_1,\ldots,w_r)+\mathbb{R}_+\{-e_1,\ldots,-e_n\})=
\{z\vert\, z\geq 0;\langle z,\ell_i\rangle\leq 1\, \forall i\}.
\end{equation}
Therefore by Eq.~(\ref{napkin-equation}) together with
Eqs.~(\ref{march6-08-3}) 
and  (\ref{march6-08-2}) we conclude that $T(P)$ is equal to 
${\rm conv}(w_1,\ldots,w_r)$, as required.
\end{proof}

If $v_1,\ldots,v_q$ are 
$\{0,1\}$-vectors, then the equality of Lemma~\ref{napkin-lemma} 
follows directly from \cite[Theorem~8]{fulkerson}; see 
also \cite{fulkerson-jctb}.

\begin{definition}\label{irpx<=1}\rm Let $A$ be a matrix with entries in
$\mathbb{N}$. The system $x\geq 0; xA\leq\mathbf{1}$ 
has the {\it integer rounding property\/} if 
$$\lceil{\rm min}\{\langle y,{\mathbf 1}\rangle \vert\, 
y\geq 0;\, Ay\geq a \}\rceil
={\rm min}\{\langle y,{\mathbf 1}\rangle \vert\, 
Ay\geq a;\, y\in\mathbb{N}^q\}
$$
for each integral vector $a$ for which ${\rm min}\{\langle
y,{\mathbf 1}\rangle \vert\,  
y\geq 0;\, Ay\geq a\}$ is finite. 
\end{definition}

If $a\in {\mathbb R}^n$, its {\it
support\/} is given by ${\rm supp}(a)=\{i|\, a_i\neq
0\}$. Note that $a=a_+-a_-$, 
where $a_+$ and $a_-$ are two non negative vectors 
with disjoint support called the {\it positive\/} and 
{\it negative\/} part of $a$ respectively. 

\begin{remark} Let $A$ be a matrix with entries in
$\mathbb{N}$. The system $x\geq 0; xA\leq\mathbf{1}$ 
has the {\it integer rounding property\/} if and only if 
$$\lceil{\rm min}\{\langle y,{\mathbf 1}\rangle \vert\, 
y\geq 0;\, Ay\geq a \}\rceil
={\rm min}\{\langle y,{\mathbf 1}\rangle \vert\, 
Ay\geq a;\, y\in\mathbb{N}^q\}
$$
for each vector $a\in\mathbb{N}^n$ for which ${\rm min}\{\langle
y,{\mathbf 1}\rangle \vert\,  y\geq 0;\, Ay\geq a\}$ is finite.
This follows decomposing an integral vector $a$ as $a=a_+-a_-$ and
noticing  
that for $y\geq 0$ we have that $Ay\geq a$ if and 
only if $Ay\geq a_+$
\end{remark}

A rational polyhedron $Q$ 
is said to have the {\it integer decomposition property\/} if for each
natural number $k$ and for each integer vector $a$ in $kQ$, $a$ is
the sum of 
$k$ integer vectors in $Q$; see \cite[pp.~66--82]{Schr2}. Recall that
$kQ$ is equal to $\{ka\vert\, a\in Q\}$. 

The next criterion will be used to describe the integer rounding 
property of the system $x\geq 0; xA\leq\mathbf{1}$ in terms 
of the normality of a certain subring.

\begin{theorem}{\rm(\cite{baum-trotter},
\cite[p.~82]{Schr2})}\label{baum-trotter-r} 
Let $A$ be a non-negative integer matrix and let $P=\{x\vert\, x\geq
0;\, xA\leq\mathbf{1}\}$. The system $xA\leq\mathbf{1};
x\geq 0$ has the integer rounding property if and only if $T(P)$ has
the integer 
decomposition property and all maximal 
integer vectors of $T(P)$ are columns of $A$ 
$($maximal with respect to $\leq$$)$. 
\end{theorem}

The next result was shown in \cite{roundp} when $A$ is the incidence
matrix of a clutter. Its proof is similar to that of \cite{roundp},
but it requires some adjustments. 

\begin{theorem}\label{round-up-char} Let $A$ be a matrix with entries in
$\mathbb{N}$ and let $v_1,\ldots,v_q$ be the columns of $A$. 
If $w_1,\ldots,w_r$ is the set of all
$\alpha\in\mathbb{N}^n$ such that $\alpha\leq v_i$ for some $i$, 
then the system $x\geq 0;\, xA\leq \mathbf{1}$ has the integer
rounding 
property if and only if the subring
$K[x^{w_1}t,\ldots,x^{w_r}t]$ is normal. 
\end{theorem}

\begin{proof} Let $P=\{x\vert\, x\geq 0;\, xA\leq\mathbf{1}\}$ and
let $T(P)$ 
be its antiblocking polyhedron. By Lemma~\ref{napkin-lemma} one has
\begin{equation}\label{key-equation}
T(P)={\rm conv}(w_1,\ldots,w_r).
\end{equation}

Let $\overline{S}$ be the integral closure of
$S=K[x^{w_1}t,\ldots,x^{w_r}t]$ 
in its field of
fractions. By the description 
of $\overline{S}$ given in Eq.~(\ref{intclos-desc-08}) one has 
$$
\overline{S}=K[\{x^at^b\, \vert\,
(a,b)\in\mathbb{Z}\mathcal{B}\cap\mathbb{R}_+\mathcal{B}\}],
$$
where $\mathcal{B}=\{(w_1,1),\ldots,(w_r,1)\}$. By
Theorem~\ref{baum-trotter-r} it suffices to prove that $S$ is normal 
if and only if $T(P)$ has the integer
decomposition property and all maximal 
integer vectors of $T(P)$ are columns of $A$ 
$($maximal with respect to $\leq$$)$. 

Assume that $S$ is normal, i.e., $S=\overline{S}$. Let $b$ be
a natural number and let $a$ be
an integer vector in $bT(P)$. Then using Eq.~(\ref{key-equation}) it is
seen that $(a,b)$ is in
$\mathbb{R}_+\mathcal{B}$. Since $S$ is normal we have
$\mathbb{R}_+\mathcal{B}\cap\mathbb{Z}\mathcal{B}=
\mathbb{N}\mathcal{B}$. In our situation one has
$\mathbb{Z}\mathcal{B}=\mathbb{Z}^{n+1}$. Hence
$(a,b)\in\mathbb{N}\mathcal{B}$  
and $a$ is the sum
of $b$ integer vectors in $T(P)$. Thus $T(P)$ has the integer
decomposition property. Assume that $a$ is a maximal integer
vector of $T(P)$. It is not hard to see that $(a,1)$ is 
in $\mathbb{R}_+\mathcal{B}$, i.e., $x^a t\in \overline{S}=S$.
Thus $(a,1)$ is a linear combination of vectors in $\mathcal{B}$ with
coefficients in 
$\mathbb{N}$. Hence $(a,1)$ is equal to $(w_j,1)$ for some $j$. There
exists $v_i$ such that $a=w_j\leq v_i$. Therefore 
by the
maximality of $a$, we get $a=v_i$ for some $i$. Thus 
$a$ is a column of $A$ as required.

Conversely assume that $T(P)$
has the integer  
decomposition property and that all maximal 
integer vectors of $T(P)$ are columns of $A$. 
Let $x^at^b\in\overline{S}$. Then $(a,b)$ is in 
the cone $\mathbb{R}_+\mathcal{B}$. Hence, using
Eq.~(\ref{key-equation}), we get $a\in
bT(P)$. Thus $a=\alpha_1+\cdots+\alpha_b$, where $\alpha_i$ is an
integral vector of $T(P)$ for all $i$. Since each $\alpha_i$ is less
than or equal to a maximal integer vector of $T(P)$, we get that 
$\alpha_i\in \{w_1,\ldots,w_r\}$. Then $x^at^b\in S$. 
This proves that $S=\overline{S}$. \end{proof}

Let $A$ be a matrix with entries in $\mathbb{N}$. Next we study the
integer rounding property of the system $x\geq 0$;
$xA\geq\mathbf{1}$. The aim is 
to establish a duality with other systems of linear inequalities.

\begin{definition}\label{irpx>=1}\rm The system $x\geq 0;
xA\geq\mathbf{1}$  
has the {\it integer rounding property\/} if 
\begin{equation}\label{irpx>=1-eq}
{\rm max}\{\langle y,{\mathbf 1}\rangle \vert\, 
y\geq 0; Ay\leq a; y\in\mathbb{N}^q\} 
=\lfloor{\rm max}\{\langle y,{\mathbf 1}\rangle \vert\, y\geq 0;
Ay\leq a\}\rfloor
\end{equation}
for each integral vector $a$ for which the right 
hand side is finite. 
\end{definition}

For any rational polyhedron $Q$ in $\mathbb{R}^n$, define its {\it
blocking polyhedron\/} $B(Q)$ by:
$$
B(Q):=\{z\in\mathbb{R}^n\vert\, z\geq 0;\, \langle z,x\rangle\geq
1\mbox{ for all 
}x\mbox{ in }Q \}.
$$
For any matrix $A$ with entries in $\mathbb{N}$, its {\it covering
polyhedron\/} 
$Q(A)$ is defined by: 
$$
Q(A):=\{x\vert\, x\geq 0; xA\geq \mathbf{1}\}.
$$
If $A$ is the incidence matrix of a clutter $\mathcal{C}$, then the 
integral vectors of $Q(A)$ correspond to vertex covers of
$\mathcal{C}$ and the 
integral vertices of $Q(A)$ are in one to one correspondence with the
minimal vertex 
covers of $\mathcal{C}$ \cite[Corollary 2.3]{reesclu}.

The blocking polyhedron of $Q(A)$ can be expressed as follows. 

\begin{lemma}\label{blocking-poly-lemma} If $Q=Q(A)$, then
$B(Q)=\mathbb{R}_+^n+{\rm conv}(v_1,\ldots,v_q)$. 
\end{lemma}
\begin{proof}
The right hand side is clearly contained in
the left hand side. Conversely take $z$ in $B(Q)$, then 
$\langle z,x\rangle\geq 1$ for all $x\in Q$ and $z\geq 0$. Let
$\ell_1,\ldots,\ell_r$ be the vertex set of $Q$. In particular 
$\langle z,\ell_i\rangle\geq 1$ for all $i$. Then 
$\langle(z,1),(\ell_i,-1)\rangle\geq 0$ for all $i$. From
\cite[Theorem 3.2]{clutters}
we get that $(z,1)$ belongs to the cone generated by  
$$
\mathcal{A}'=\{e_1,\ldots,e_n,(v_1,1),\ldots,(v_q,1)\}.
$$ 
Thus $z$ is
in $\mathbb{R}_+^n+{\rm conv}(v_1,\ldots,v_q)$. This completes the
proof of the asserted equality.
\end{proof}

The next criterion complements Theorem~\ref{baum-trotter-r}.

\begin{theorem}{\rm(\cite{baum-trotter},
\cite[p.~82]{Schr2})}\label{baum-trotter-r-dual}\  
The system $x\geq 0; xA\geq\mathbf{1}$ has the integer rounding
property if and only if  
the blocking polyhedron $B(Q)$ of $Q=Q(A)$ has the integer
decomposition property and all minimal
integer vectors of $B(Q)$ are columns of $A$ 
$($minimal with respect to $\leq$$)$. 
\end{theorem}

Recall that a set 
$\mathcal{A}\subset\mathbb{Z}^n$ is called a {\it Hilbert basis\/} if
$\mathbb{N}\mathcal{A}=\mathbb{R}_+\mathcal{A}\cap\mathbb{Z}^n$. Note
that if $\mathcal{A}$ is a Hilbert basis, then the semigroup ring 
$K[\mathbb{N}\mathcal{A}]$ is normal.

\begin{proposition}\label{jun25-08} Let $I=(x^{v_1},\ldots,x^{v_q})$
be a monomial 
ideal and let $v_i^*=\mathbf{1}-v_i$. Then $R[It]$ is normal if and
only if the set  
$$
\Gamma=\{-e_1,\ldots,-e_n,(v_1^*,1),\ldots,(v_q^*,1)\}
$$
is a Hilbert basis.
\end{proposition}

\begin{proof} Let
$\mathcal{A}'=\{e_1,\ldots,e_n,(v_1,1),\ldots,(v_q,1)\}$. Assume that
$R[It]$ is normal. Then $\mathcal{A}'$ is a Hilbert basis. Let 
$(a,b)$ be an integral vector in $\mathbb{R}_+\Gamma$, 
with $a\in\mathbb{Z}^n$ and $b\in\mathbb{Z}$. Then we can write
$$ 
(a,b)=\mu_1(-e_1)+\cdots+\mu_n(-e_n)+\lambda_1(v_1^*,1)+\cdots+
\lambda_q(v_q^*,1), 
$$
where $\mu_i\geq 0$ and $\lambda_j\geq 0$ for all $i,j$. Therefore
$$
-(a,b)+b\mathbf{1}=\mu_1e_1+\cdots+\mu_ne_n+\lambda_1(v_1,-1)+\cdots+
\lambda_q(v_q,-1),
$$
where $\mathbf{1}=e_1+\cdots+e_n$. This equality is equivalent to
$$
-(a,-b)+b\mathbf{1}=\mu_1e_1+\cdots+\mu_ne_n+\lambda_1(v_1,1)+\cdots+
\lambda_q(v_q,1).
$$
As $\mathcal{A}'$ is a Hilbert basis we can write
$$
-(a,-b)+b\mathbf{1}=\mu_1'e_1+\cdots+\mu_n'e_n+\lambda_1'(v_1,1)+\cdots+
\lambda_q'(v_q,1),
$$
where $\mu_i'\in\mathbb{N}$ and $\lambda_j'\in\mathbb{N}$ for all $i,j$. Thus
$(a,b)\in\mathbb{N}\Gamma$. This proves that $\Gamma$ is a Hilbert
basis. The converse can be shown using similar arguments.
\end{proof}

A {\it clutter\/} $\mathcal{C}$ with finite vertex set
$X=\{x_1,\ldots,x_n\}$ is a 
family of subsets of $X$, 
called edges, none of which is included in another. Let
$f_1,\ldots,f_q$ be 
the edges of $\mathcal{C}$ and let $v_k=\sum_{x_i\in f_k}e_i$ be the
{\it characteristic vector\/} of $f_k$. The {\it incidence matrix\/} 
of $\mathcal{C}$ is the $n\times q$ matrix with 
column vectors $v_1,\ldots,v_q$.

\begin{definition} Let $A=(a_{ij})$ be a matrix with entries in
$\{0,1\}$. Its {\it dual\/} 
is the matrix $A^*=(a_{ij}^*)$, where $a_{ij}^*=1-a_{ij}$.
\end{definition}

The following duality is valid for incidence matrices of clutters. It
will be used later to establish a 
duality theorem for monomial subrings.

\begin{theorem}\label{duality-irp-clutters} Let $A$ be the incidence
matrix of a clutter and let $v_1,\ldots,v_q$ be its column vectors. 
If $v_i^*=\mathbf{1}-v_i$ and $A^*$ is the matrix with column vectors
$v_1^*,\ldots,v_q^*$, then the system $x\geq 0;xA\geq\mathbf{1}$
has the integer 
rounding property if and 
only if the system $x\geq 0;xA^*\leq\mathbf{1}$ has the integer
rounding property.
\end{theorem}

\begin{proof} Consider $Q=\{x\vert x\geq 0;xA\geq \mathbf{1}\}$ and 
$P^*=\{x\vert x\geq 0;xA^*\leq\mathbf{1}\}$. Let $w_1^*,\ldots,w_s^*$
be the set of all $\alpha\in\mathbb{N}^n$ such that $\alpha\leq
v_i^*$ for some $i$. Then, using Lemmas~\ref{blocking-poly-lemma} and
\ref{napkin-lemma}, we 
obtain that the blocking polyhedron of $Q$ and the antiblocking 
polyhedron of $P^*$ are given by 
$$
B(Q)=\mathbb{R}_+^n+{\rm conv}(v_1,\ldots,v_q)\ \mbox{ and }\  
T(P^*)={\rm conv}(w_1^*,\ldots,w_s^*)
$$
respectively.

$\Rightarrow$) By Theorem~\ref{baum-trotter-r} it suffices to show
that $T(P^*)$ 
has the integer decomposition property and all maximal 
integer vectors of $T(P^*)$ are columns of $A^*$. Let $b$ be an
integer and let $a$ be an integer vector in $bT(P^*)$. Then we
can write
$$
a=b(\lambda_1w_1^*+\cdots+\lambda_sw_s^*),\ \ \
(\textstyle\sum_i\lambda_i=1;\lambda_i\geq 0).
$$
For each $1\leq i\leq s$ there is $v_{j_i}^*$ in
$\{v_1^*,\ldots,v_q^*\}$ such that $w_i^*\leq
v_{j_i}^*$. Thus for each $i$ we can write 
$\mathbf{1}-w_i^*=v_{j_i}+\delta_i$, where
$\delta_i\in\mathbb{N}^n$. Therefore 
$$
\mathbf{1}-a/b=
\lambda_1(v_{j_1}+\delta_1)+\cdots+\lambda_s(v_{j_s}+\delta_s).
$$
This means that $\mathbf{1}-a/b\in B(Q)$, i.e., $b\mathbf{1}-a$ is an
integer vector in $bB(Q)$. Hence by Theorem~\ref{baum-trotter-r-dual}
we can write $b\mathbf{1}-a=\alpha_1+\cdots+\alpha_b$ for some
$\alpha_1,\ldots,\alpha_b$ integer vectors in $B(Q)$, and for each
$\alpha_i$ there is $v_{k_i}$ in $\{v_1,\ldots,v_q\}$ such that
$v_{k_i}\leq\alpha_i$. Thus $\alpha_i=v_{k_i}+\epsilon_i$ for some
$\epsilon_i\in\mathbb{N}^n$ and consequently:
$$
a=(\mathbf{1}-v_{k_1})+\cdots+(\mathbf{1}-v_{k_b})-c
=v_{k_1}^*+\cdots+v_{k_b}^*-c,
$$ 
where $c=(c_1,\ldots,c_n)\in\mathbb{N}^n$. Notice that
$v_{k_1}^*+\cdots+v_{k_b}^*\geq c$ because $a\geq 0$. If $c_1\geq 1$,
then the first entry of $v_{k_i}^*$ 
is non-zero for some $i$ and we can write
$$
a=v_{k_1}^*+\cdots+v_{k_{i-1}}^*+(v_{k_i}^*-e_1)+v_{k_{i+1}}^*+
\cdots+v_{k_b}^*-(c-e_1)
$$
Since $v_{k_i}^*-e_1$ is again in $\{w_1^*,\ldots,w_s^*\}$, we can
apply this argument recursively to obtain that $a$ is the sum of 
$b$ integer vectors in $\{w_1^*,\ldots,w_s^*\}$. This proves that
$T(P^*)$ has the integer decomposition property. Let $a$ be a maximal
integer vector of $T(P^*)$. Since the vectors $w_1^*,\ldots,w_s^*$
have entries in $\{0,1\}$, we get
$T(P^*)\cap\mathbb{Z}^n=\{w_1^*,\ldots,w_s^*\}$. Then $a=w_i^*$ for
some $i$. As $w_i^*\leq v_j^*$ for some $j$, we conclude that
$a=v_j^*$, i.e., $a$ is a column of $A^*$, as required. 

$\Leftarrow$) According to \cite{poset} the system $x\geq 0;
xA\geq\mathbf{1}$ has the integer rounding property if and only if 
$R[It]$ is normal. Thus by Proposition~\ref{jun25-08} we need only
show that the set
$\Gamma=\{-e_1,\ldots,-e_n,(v_1^*,1),\ldots,(v_q^*,1)\}$ is a Hilbert
basis. Let  
$(a,b)$ be an integral vector in $\mathbb{R}_+\Gamma$, 
with $a\in\mathbb{Z}^n$ and $b\in\mathbb{Z}$. Then we can write
$$ 
(a,b)=\mu_1(-e_1)+\cdots+\mu_n(-e_n)+\lambda_1(v_1^*,1)+\cdots+
\lambda_q(v_q^*,1),
$$
where $\mu_i\geq 0$, $\lambda_j\geq 0$ for all $i,j$. Hence
$A^*\lambda\geq a$, where $\lambda=(\lambda_i)$. By hypothesis 
the system $x\geq 0;xA^*\leq\mathbf{1}$ has the integer rounding
property. Then one has
$$
b\geq \lceil{\rm min}\{\langle y,{\mathbf 1}\rangle \vert\, 
y\geq 0;\, A^*y\geq a \}\rceil
={\rm min}\{\langle y,{\mathbf 1}\rangle \vert\, 
A^*y\geq a;\, y\in\mathbb{N}^q\}=\langle{y_0},\mathbf{1}\rangle
$$
for some $y_0=(y_i)\in\mathbb{N}^q$ such that
$|y_0|=\langle{y_0},\mathbf{1}\rangle\leq b$ and $a\leq A^*y_0$. Then
$$
a=y_1v_1^*+\cdots+y_qv_q^*-\delta_1e_1-\cdots-\delta_ne_n,
$$
where $\delta_1,\ldots,\delta_n$ are in $\mathbb{N}$. Hence we can
write
$$ 
(a,b)=y_1(v_1^*,1)+\cdots+y_{q-1}(v_{q-1}^*,1)+(y_q+b-|y_0|)(v_q^*,1)
-(b-|y_0|)v_q^*-\delta,
$$
where $\delta=(\delta_i)$. As the entries of $A^*$ are in $\mathbb{N}$,
the vector $-v_q^*$ can be written as a non-negative integer
combination of $-e_1,\ldots,-e_n$. Thus
$(a,b)\in\mathbb{N}\Gamma$. This proves that $\Gamma$ is a
Hilbert basis. 
\end{proof}

We come to one of the main result of this section. It establishes a 
duality for monomial subrings.

\begin{theorem}\label{duality-subrings-irp} Let $A$ be the incidence
matrix of a clutter, let $v_1,\ldots,v_q$ be its column vectors and
let $v_i^*=\mathbf{1}-v_i$. If 
$w_1^*,\ldots,w_s^*$ is the set of all $\alpha\in\mathbb{N}^n$ such 
that $\alpha\leq v_i^*$ for some $i$, then the following conditions
are equivalent: 
\begin{itemize}
\item[(a)] $R[It]$ is normal, where $I=(x^{v_1},\ldots,x^{v_q})$. 
\item[(b)] $S^*=K[x^{w_1^*}t,\ldots,x^{w_s^*}t]$ is normal. 
\item[(c)] $\{-e_1,\ldots,-e_n,(v_1^*,1),\ldots,(v_q^*,1)\}$ is a
Hilbert basis. 
\item[(d)] $x\geq 0;xA\geq\mathbf{1}$ has the integer rounding
property.
\item[(e)] $x\geq 0;xA^*\leq\mathbf{1}$ has the integer rounding
property.
\end{itemize}
\end{theorem}

\begin{proof}(a) $\Leftrightarrow$ (c): This was shown in
Proposition~\ref{jun25-08}. (a) $\Leftrightarrow$ (d): This is one of
the main results of \cite{poset} and is valid for arbitrary monomial
ideals. (b) $\Leftrightarrow$ (e): This was
shown in Theorem~\ref{round-up-char}. (d) $\Leftrightarrow$ (e): 
This follows from Theorem~\ref{duality-irp-clutters}.
\end{proof}

To illustrate the usefulness of this duality, below we show various 
results that follow from there.

\begin{definition} Let $\mathcal{C}$ be a clutter on the vertex set
$X=\{x_1,\ldots,x_n\}$. The {\it edge ideal\/} of $\mathcal{C}$,
denoted by $I(\mathcal{C})$, is the ideal of $R$ generated by 
all monomials $x_e=\prod_{x_i\in e}x_i$ such that $e$ is an edge of
$\mathcal{C}$. The {\it dual\/} $I^*$ of an edge ideal $I$ is the
ideal of $R$ 
generated by all monomials $x_1\cdots x_n/x_e$ such that $e$ is an
edge of $\mathcal{C}$.
\end{definition}

\begin{corollary}{\rm (\cite[Theorem~2.10]{perfect})}\label{jun26-08}
Let $\mathcal{C}$ be a clutter and let $A$ be its incidence matrix. If
$P=\{x\vert x\geq 0;xA\leq\mathbf{1}\}$ is an 
integral polytope and $I=I(\mathcal{C})$, then
\begin{itemize}
\item[(i)] $R[I^*t]$ is normal. 
\item[(ii)] $S=K[x^{w_1}t,\ldots,x^{w_r}t]$ is normal.
\end{itemize}
\end{corollary}

\begin{proof} Since $P$ has only integral vertices, by a result 
of Lov\'asz \cite{lovasz}
the system $x\geq 0;xA\leq\mathbf{1}$ is totally dual
integral, i.e., the minimum in 
the LP-duality equation
\begin{equation}\label{jun6-2-03}
{\rm max}\{\langle a,x\rangle \vert\, x\geq 0; xA\leq \mathbf{1}\}=
{\rm min}\{\langle y,\mathbf{1}\rangle \vert\, y\geq 0; Ay\geq a\} 
\end{equation}
has an integral optimum solution $y$ for each integral vector $a$ with 
finite minimum. In particular the system $x\geq 0;xA\leq\mathbf{1}$
satisfies the integer rounding property. Therefore $R[I^*t]$ and 
$K[x^{w_1}t,\ldots,x^{w_r}t]$ are normal by
Theorem~\ref{duality-subrings-irp}.  
\end{proof}

This result is related to the theory of perfect graphs. Indeed if $P$ is
integral, the $w_i's$ correspond
to the cliques (complete subgraphs) of a perfect graph $H$
\cite{chvatal,lovasz}, 
and the $v_i^*$'s correspond to the minimal vertex covers of the
complement of $H$. The normality assertion of part (ii) is well known
and it can also be 
shown directly using the fact that the system $x\geq
0;xA\leq\mathbf{1}$ is TDI if $P$ 
is integral, where TDI stands for Totally Dual Integral 
(see \cite{Schr2}).

\begin{corollary}\label{jun29-08} Let $B_1,\ldots,B_q$ be the
collection of basis of a matroid $M$ with vertex set $X$ and let
$v_1,\ldots,v_q$ be their characteristic vectors. If $A$ is the matrix
with column vectors $v_1,\ldots,v_q$, then all systems
$$
x\geq 0;xA\geq\mathbf{1},\ \ \ x\geq 0;xA^*\geq\mathbf{1},\ \ \ 
x\geq 0;xA\leq\mathbf{1},\ \ \ x\geq 0;xA^*\leq\mathbf{1}
$$
have the integer rounding property.
\end{corollary}
\begin{proof} Consider the basis monomial ideal
$I=(x^{v_1},\ldots,x^{v_q})$ of the matroid $M$. By
\cite[Theorem~2.1.1]{oxley}, the collection of
basis of the dual matroid $M^*$ of $M$ is given by 
$X\setminus B_1,\ldots,X\setminus B_q$. Now, the basis monomial
ideal of a matroid is normal \cite[Corollary~3.8]{matrof}, thus the
result follows at once from the duality given in
Theorem~\ref{duality-subrings-irp}.  
\end{proof}

\begin{corollary}\label{forgraphs-In-iff-dualIn} Let $G$ be a
connected graph and let $I=I(G)$ be its edge
ideal. Then $R[It]$ is normal if and only if $R[I^*t]$ is normal.
\end{corollary}

\begin{proof} By \cite{roundp} the system $x\geq 0;xA\geq\mathbf{1}$
has the integer rounding property if and only if the system 
$x\geq 0;xA\leq\mathbf{1}$ does. Therefore the result follows at once
using Theorem~\ref{duality-subrings-irp}.
\end{proof}

This result is valid even if the graph is not connected but its proof
requires to use the fact that $R[It]$ is normal if and only 
if the extended Rees algebra $R[It,t^{-1}]$ is normal and the fact
that $R[It,t^{-1}]$ is isomorphic to $S=K[x^{w_1},\ldots,x^{w_r}]$
when $I$ is the edge ideal of a graph (see \cite{roundp}). 
The next example shows that Corollary~\ref{forgraphs-In-iff-dualIn} 
 does not extends to arbitrary uniform clutters.

\begin{example}\label{forclutters-In-iff-dualIn-false} 
Consider the clutter $\mathcal{C}$ whose incidence
matrix $A$ is the transpose of the matrix:
\begin{small}
$$
\left[
\begin{array}{cccccccccc}
0& 0& 1& 1& 0& 1& 1& 1& 1& 1\\ 
0& 0& 1& 0& 1& 1& 1& 1& 1& 1\\  
0& 1& 1& 0& 0& 1& 1& 1& 1& 1\\  
1& 1& 0& 0& 0& 1& 1& 1& 1& 1\\  
0& 1& 1& 0& 1& 0& 1& 1& 1& 1\\ 
1& 1& 1& 1& 1& 0& 0& 1& 1& 0\\ 
1& 1& 1& 1& 1& 0& 0& 1& 0& 1\\ 
1& 1& 1& 1& 1& 0& 1& 1& 0& 0\\  
1& 1& 1& 1& 1& 1& 1& 0& 0& 0\\ 
1& 1& 1& 1& 0& 0& 1& 1& 0& 1
\end{array}
\right]
$$
\end{small}
Let $I=I(\mathcal{C})$ be the edge ideal of $\mathcal{C}$. Note that
all edges of $\mathcal{C}$ have $7$ vertices.  
Using {\it Normaliz\/} \cite{B} it is seen that $R[It]$ is normal and
that $R[I^*t]$ is not normal. 
\end{example}

Let $\mathcal{C}$ be a clutter with vertex set $X$. A vertex $x$ of
$\mathcal{C}$ is called {\it isolated\/} if $x$ does not occur in any
edge of $\mathcal{C}$. A subset $C\subset X$ is a 
{\it minimal vertex cover\/} of $\mathcal{C}$ if: 
($\mathrm{c}_1$) every edge of $\mathcal{C}$ contains at least one
vertex of $C$,  
and ($\mathrm{c}_2$) there is no proper subset of $C$ with the first 
property. If $C$ only satisfies condition ($\mathrm{c}_1$), then $C$ is 
called a {\it vertex cover\/} of $\mathcal{C}$. The {\it Alexander
dual\/} of 
$\mathcal{C}$, denoted by $\mathcal{C}^\vee$, is the clutter 
whose edges are the minimal vertex covers of $\mathcal{C}$. The edge
ideal of $\mathcal{C}^\vee$, denoted by $I(\mathcal{C})^\vee$, 
is called the {\it Alexander dual\/} of $I(\mathcal{C})$. In
combinatorial optimization the 
Alexander dual of
a clutter is referred to as the {\it blocker\/} of the clutter
\cite{Schr2}. 

\begin{proposition}\label{el-vegetariano} Let $G$ be a graph without
isolated vertices and let $G'$ be its complement. Then 
$I(G')^\vee=I(G)^*$ if and only if $G$ is triangle free.
\end{proposition}

\begin{proof} $\Rightarrow$) Let $X=\{x_1,\ldots,x_n\}$ be the vertex
set of $G$. Assume that $G$  has a triangle
$\mathcal{C}_3=\{x_1,x_2,x_3\}$, i.e., $\{x_i,x_j\}$ are edges of $G$
for $1\leq i<j\leq 3$. Clearly we may assume $n\geq 4$. Notice that
$C'=\{x_4,\ldots,x_n\}$ is a vertex cover of $G'$, i.e., $x_4\cdots
x_n$ belongs to $I(G')^\vee$ and consequently it belongs to $I(G)^*$,
a contradiction because $I(G)^*$ is generated by monomials of degree
$n-2$.

$\Leftarrow$) Let $x^a=x_1\cdots x_r$ be a minimal generator of
$I(G')^\vee$. Then $C=\{x_1,\ldots,x_r\}$ is a minimal vertex cover 
of $G'$. Hence $X\setminus C$ is a maximal complete subgraph of $G$.
Thus by hypothesis $X\setminus C$ is an edge of $G$, i.e., $x^a\in
I(G)^*$. This proves the inclusion $I(G')^\vee\subset I(G)^*$.
Conversely, let $x^a$ be a minimal generator of $I(G)^*$. There is an
edge $\{x_1,x_2\}$ of $G$ such that $x^a=x_3\cdots x_n$. Every edge of
$G'$ must intersect $C=\{x_3,\ldots,x_n\}$, i.e., $x^a\in I(G')^\vee$.
\end{proof}

This formula applies for instance if $G$ is a bipartite graph.

\begin{corollary}\label{july4-08} Let $G$ be a free triangle graph
without isolated 
vertices. Then $R[I(G)t]$ is normal if and only if $R[I(G')^\vee t]$
is normal. 
\end{corollary}

\begin{proof} It follows directly from
Corollary~\ref{forgraphs-In-iff-dualIn} and
Proposition~\ref{el-vegetariano}. 
\end{proof}

In \cite{perfect} it is shown that the Alexander dual of the edge
ideal of a perfect graph is always normal (cf.
Corollary~\ref{jun26-08}(i)). 
To the best of our knowledge the following is the first example of an
edge ideal of a graph whose Alexander dual is not normal.

\begin{example}\label{contraec} Let $G$ be the graph consisting of
two vertex disjoint 
odd cycles of length $5$ and let $G'$ be its complement. According to
\cite{bowtie} the Rees algebra of $I(G)$ is not normal. Thus
$R[I(G')^\vee{t}]$ is not normal by Corollary~\ref{july4-08}.
\end{example}
\begin{definition}\label{mfmc-def}\rm 
A clutter $\mathcal C$ satisfies the {\it max-flow min-cut\/}
(MFMC) 
property if both sides 
of the LP-duality equation
\begin{equation}\label{jun6-2-03-1}
{\rm min}\{\langle a,x\rangle \vert\, x\geq 0; xA\geq{\mathbf 1}\}=
{\rm max}\{\langle y,{\mathbf 1}\rangle \vert\, y\geq 0; Ay\leq a\} 
\end{equation}
have integral optimum solutions $x$ and $y$ for each non-negative 
integral vector $a$. 
\end{definition}

\begin{corollary}{\rm(\cite[Theorem~{3.4}]{clutters})}\label{jun29-08-1}
Let $A$ be the incidence matrix of
a clutter $\mathcal{C}$ and let $I=I(\mathcal{C})$ be its edge ideal.
Then $\mathcal{C}$ satisfies the max-flow min-cut
property if and only if $Q(A)$ is integral and $R[It]$ is normal.
\end{corollary}

\begin{proof} Notice that if $\mathcal{C}$ has the max-flow min-cut
property, then $Q(A)$ is integral \cite[Corollary~22.1c]{Schr}.
Therefore the result follows directly from Eqs.~(\ref{irpx>=1-eq}),
(\ref{jun6-2-03-1}), and Theorem~\ref{duality-subrings-irp}.
\end{proof}

We now turn our attention to the integer rounding 
property of systems of the form $xA\leq\mathbf{1}$. 

\begin{definition}\label{irp}\rm Let $A$ be a matrix with entries in
$\mathbb{N}$. The system $xA\leq\mathbf{1}$ is said to have 
the {\it integer rounding property\/} if 
$$\lceil{\rm min}\{\langle y,{\mathbf 1}\rangle \vert\, 
y\geq 0;\, Ay=a \}\rceil
={\rm min}\{\langle y,{\mathbf 1}\rangle \vert\, 
Ay=a;\, y\in\mathbb{N}^q\}
$$
for each integral vector $a$ for which ${\rm min}\{\langle
y,{\mathbf 1}\rangle \vert\,  
y\geq 0;\, Ay=a\}$ is finite. 
\end{definition}

The next result is just a reinterpretation of an unpublished result
of Giles and 
Orlin \cite[Theorem~22.18]{Schr} that characterizes the integer
rounding property in terms of Hilbert bases.

\begin{proposition}\label{giles-orlin} Let $v_1,\ldots,v_q$ be the
column vectors of a non-negative integer matrix $A$ and let
$A(P)$ be the Ehrhart ring of 
$P={\rm conv}(0,v_1,\ldots,v_q)$.
Then the system $xA\leq\mathbf{1}$ has the integer rounding 
property if and only if 
$$
K[x^{v_1}t,\ldots,x^{v_q}t,t]=A(P).
$$
\end{proposition}

\begin{proof} By \cite[Theorem~22.18]{Schr}, we have that the system
$xA\leq\mathbf{1}$ has the integer rounding  
property if and only if the set 
$\mathcal{B}=\{(v_1,1),\ldots,(v_q,1),(0,1)\}$ is
a Hilbert basis. Thus the proposition follows readily 
by noticing the equality
$$A(P)=K[\{x^at^b\vert(a,b)\in\mathbb{R}_+\mathcal{B}\cap\mathbb{Z}^{n+1}\}]$$
and the inclusion $K[x^{v_1}t,\ldots,x^{v_q}t,t]\subset A(P)$. \end{proof}

\begin{theorem}\label{jun5-08} Let $\mathcal{A}=\{v_1,\ldots,v_q\}$
be the set of  
column vectors of a matrix $A$ with entries in $\mathbb{N}$. 
If the system
$xA\leq\mathbf{1}$ has the integer rounding property, then 
\begin{itemize}
\item[(a)] $K[F]$ is normal, where $F=\{x^{v_1},\ldots,x^{v_q} \}$, 
and \item[(b)] $\mathbb{Z}^n/\mathbb{Z}\mathcal{A}$ is a torsion-free group. 
\end{itemize}
The converse holds if $|v_i|=d$ for all $i$.
\end{theorem}

\begin{proof} For use below we set
$\mathcal{B}=\{(v_1,1),\ldots,(v_q,1),(0,1)\}$. First we prove (a).
Let $x^a\in \overline{K[F]}$. Then $a\in\mathbb{Z}\mathcal{A}$ and we
can write  
$$
a=\lambda_1v_1+\cdots+\lambda_qv_q,
$$
for some $\lambda_1,\ldots,\lambda_q$ in $\mathbb{R}_+$. Hence 
$$
(a,\textstyle\lceil\sum_i\lambda_i
\rceil)=\lambda_1(v_1,1)+\cdots+\lambda_q(v_q,1)+
\delta(0,1),
$$
where $\delta\geq 0$. Therefore by Proposition~\ref{giles-orlin},
there are $\lambda_1',\ldots\lambda_q'\in\mathbb{N}$ and
$\delta'\in\mathbb{N}$ such that
$$
(a,\textstyle\lceil\sum_i\lambda_i
\rceil)=\lambda_1'(v_1,1)+\cdots+\lambda_q'(v_q,1)+
\delta'(0,1),
$$
Thus $x^a\in K[F]$, as required. Next we show (b). From
Proposition~\ref{giles-orlin}, we get
$$
K[x^{v_1}t,\ldots,x^{v_q}t,t]=A(P).
$$
Hence using \cite[Theorem~3.9]{ehrhart} we obtain that the group
$M=\mathbb{Z}^{n+1}/\mathbb{Z}\mathcal{B}$ is torsion free. Let
$\overline{a}$ be an element of 
$T(\mathbb{Z}^n/\mathbb{Z}\mathcal{A})$, the torsion subgroup of
$\mathbb{Z}^n/\mathbb{Z}\mathcal{A}$. Thus there is a positive
integer $s$ so that  
$$
sa=\lambda_1v_1+\cdots+\lambda_qv_q
$$
for some $\lambda_1,\ldots,\lambda_q$ in $\mathbb{Z}$. From the
equality 
$$
s(a,|a|)=\lambda_1(v_1,1)+\cdots+\lambda_q(v_q,1)+(s|a|-
\lambda_1-\cdots-\lambda_q)(0,1)
$$
we obtain that the image of $(a,|a|)$ in $M$, denoted by
$\overline{(a,|a|)}$,  is 
a torsion element, i.e., $\overline{(a,|a|)}\in T(M)=(\overline{0})$.
Hence it is readily seen that $a\in \mathbb{Z}\mathcal{A}$, i.e., 
$\overline{a}=\overline{0}$. Altogether we have
$T(\mathbb{Z}^n/\mathbb{Z}\mathcal{A})=(0)$. 

Conversely assume that $|v_i|=d$ for all $i$ and that (a) and (b)
hold. We need only show that $\mathcal{B}$ is a Hilbert basis. Let
$(a,b)$ be an integral vector in $\mathbb{R}_+\mathcal{B}$, where
$a\in\mathbb{N}^n$ and $b\in \mathbb{N}$. Then we can write
\begin{equation}\label{april8-08}
(a,b)=\lambda_1(v_1,1)+\cdots+\lambda_q(v_q,1)+\mu(0,1),
\end{equation}
for some $\lambda_1,\ldots,\lambda_q,\mu$ in $\mathbb{Q}_+$. Hence
using this equality together with (b) gives that $a$ is in 
$\mathbb{R}_+\mathcal{A}\cap\mathbb{Z}\mathcal{A}$. Hence 
$x^a\in\overline{K[F]}=K[F]$, i.e., $a\in\mathbb{N}\mathcal{A}$. Then
we can 
write
$$
a=\eta_1v_1+\cdots+\eta_qv_q
$$
for some $\eta_1,\ldots,\eta_q$ in $\mathbb{N}^n$. Since 
$|v_i|=d$ for all $i$, one has $\sum_i\lambda_i=\sum_i\eta_i$.
Therefore using Eq.~(\ref{april8-08}), we get $\mu\in\mathbb{N}$.
Consequently from the equality
$$
(a,b)=\eta_1(v_1,1)+\cdots+\eta_q(v_q,1)+\mu(0,1),
$$
we conclude that $(a,b)\in\mathbb{N}\mathcal{B}$. This proves that 
$\mathcal{B}$ is a Hilbert basis. \end{proof}

\begin{corollary}\label{jun8-08-1} Let $A$ be the incidence matrix of a 
connected graph $G$. Then the system $xA\leq\mathbf{1}$ has the
integer rounding property if and only if $G$ is a bipartite graph.
\end{corollary}

\begin{proof} $\Rightarrow$) Let $\mathcal{A}=\{v_1,\ldots,v_q\}$ be
the set 
of columns of $A$. If $G$ is not bipartite, then according to 
\cite[Corollary~3.4]{accota} one has 
$\mathbb{Z}^n/\mathbb{Z}\mathcal{A}\simeq\mathbb{Z}_2$, a
contradiction to Theorem~\ref{jun5-08}(b). 

$\Leftarrow$) By \cite[Theorem~2.15, Corollary~3.4]{accota} we get
that the ring $K[x^{v_1},\ldots,x^{v_q}]$ is normal and that 
$\mathbb{Z}^n/\mathbb{Z}\mathcal{A}\simeq\mathbb{Z}$. Thus by
Theorem~\ref{jun5-08} the system $xA\leq\mathbf{1}$ has the integer
rounding property, as required. This part of the proof also follows
directly from 
the fact that the incidence matrix of a bipartite graph is totally
unimodular. Indeed, since $A$ is totally unimodular, both problems of
the 
LP-duality
equation 
$$
\max\{\langle a,x\rangle\vert\, xA\leq\mathbf{1}\}=\min\{\langle 
y,\mathbf{1}\rangle\vert\, y\geq 0; Ay=a\} 
$$
have integral optimum solutions for each integral vector $a$ for
which the minimum is finite, see 
\cite[Corollary~19.1a]{Schr}. Thus the system $xA\leq\mathbf{1}$ has
the integer 
rounding property. \end{proof}

\begin{corollary}\label{jun18-08} Let $A$ be the incidence matrix of
a clutter 
$\mathcal{C}$. If $\mathcal{C}$ is uniform and has the max-flow
min-cut 
property, 
then the system $xA\leq\mathbf{1}$ has the
integer rounding property.
\end{corollary}

\begin{proof} Since all edges of $\mathcal{C}$ have the same size, it
suffices to observe that conditions (a) and (b) of  
Theorem~\ref{jun5-08} are satisfied because of
\cite[Theorem~3.6]{mfmc}. \end{proof}  

\section{The canonical module and the
$a$-invariant}\label{can-mod-device}

Let $R=K[x_1,\ldots,x_n]$ be a polynomial ring over an arbitrary
field $K$ and 
let $K[F]=K[x^{v_1},\ldots,x^{v_q}]$ be a homogeneous monomial subring, i.e., 
there exists $0\neq x_0\in\mathbb{Q}^n$
satisfying $\langle x_0,v_i\rangle=1$ for all $i$. Then 
$K[F]$ is a standard graded $K$-algebra with the grading induced by
declaring that a monomial  $x^{a}\in K[F]$ has degree $i$ if
and only if $\langle a,x_0\rangle=i$. 
Recall that the $a$-{\it invariant\/} of $K[F]$, denoted by $a(K[F])$, is 
the degree as a rational 
function of the Hilbert series 
of $K[F]$, see for instance \cite[p.~99]{monalg}. Let $H$ and
$\varphi$ be the Hilbert function and the Hilbert polynomial of $K[F]$
respectively. The index of regularity of $K[F]$, denoted by ${\rm
reg}(K[F])$, is the least positive integer such that $H(i)=\varphi(i)$
for $i\geq {\rm reg}(K[F])$. The $a$-invariant plays a fundamental role
in algebra and geometry because one has: ${\rm reg}(K[F])=0$ if
$a(K[F])<0$ and ${\rm reg}(K[F])=a(K[F])+1$ otherwise
\cite[Corollary~4.1.12]{monalg}.

If $K[F]$ is 
Cohen-Macaulay and $\omega_{K[F]}$ is the canonical 
module of $K[F]$, then
\begin{equation}\label{princeton-fall-07}
a(K[F])=-{\rm min}\{\, i\, \vert\, (\omega_{K[F]})_i\neq 0\},
\end{equation}
see \cite[p.~141]{BHer} and \cite[Proposition~4.2.3]{monalg}. 
This formula applies if $K[F]$ is normal because normal 
monomial subrings are Cohen-Macaulay \cite{Ho1}.  
If $K[F]$ is normal,
then by a formula of 
Danilov and Stanley (see \cite[Theorem~6.3.5]{BHer} and \cite{Dan}) 
the canonical module of $K[F]$ is 
the ideal given by  
\begin{equation}\label{ejc8}
\omega_{K[F]}=(\{x^a\vert\,
a\in \mathbb{N}{\mathcal A}\cap ({\mathbb R}_+{\mathcal A})^{\rm o}\}),
\end{equation}
where  ${\mathcal A}=\{v_1,\ldots,v_q\}$ and $({\mathbb R}_+{\mathcal
A})^{\rm o}$ is the interior of  
${\mathbb R}_+{\mathcal A}$ relative to 
${\rm aff}({\mathbb R}_+{\mathcal A})$, 
the affine hull of ${\mathbb R}_+{\mathcal A}$. 

The {\it dual cone\/} of ${\mathbb R}_+{\mathcal A}$
is the polyhedral cone given by
$$
(\mathbb{R}_+{\mathcal A})^*=\{x\, \vert\, \langle x,y\rangle\geq
0;\, \forall\, 
y\in\mathbb{R}_+{\mathcal A}\}.
$$
A set $\mathcal{H}\subset\mathbb{R}^n\setminus\{0\}$ is called
an {\it integral 
basis\/} of $(\mathbb{R}_+{\mathcal A})^*$ if  $(\mathbb{R}_+{\mathcal
A})^*=\mathbb{R}_+{\mathcal H}$ and $\mathcal{H}\subset\mathbb{Z}^n$.
Let $0\neq a\in\mathbb{R}^n$. In what follows $H_{a}^+$ denotes the
closed halfspace  
$H_a^+=\{x\vert\, \langle
x,a\rangle\geq 0\}$ and $H_a$ stands for the hyperplane through the
origin with normal 
vector $a$.

The next result gives a general technique to compute the canonical
module and 
the $a$-invariant of a wide class of monomial subrings. Another
technique is given in \cite{shiftcon}. In
Section~\ref{section-on-canmod} we give some more precise expressions
for the canonical module and the $a$-invariant of special families of
monomial subrings arising from integer rounding properties.

\begin{theorem}\label{canmod-tdi} Let $c_1,\ldots,c_r$ be an integral
basis of $(\mathbb{R}_+{\mathcal A})^*$ and let $b=(b_i)$ be the 
$\{0,-1\}$-vector given by $b_i=0$ if $\mathbb{R}_+{\mathcal
A}\subset H_{c_i}$ and $b_i=-1$ if  
$\mathbb{R}_+{\mathcal A}\not\subset H_{c_i}$. If $\mathbb{N}{\mathcal
A}=\mathbb{Z}^n\cap \mathbb{R}_+{\mathcal A}$ and $B$ is the matrix
with column vectors $-c_1,\ldots,-c_r$, then  
\begin{itemize}
\item[(a)] $\omega_{K[F]}=(\{x^a\vert\,
a\in \mathbb{Z}^n\cap \{x\vert xB\leq b\})$.
\item[(b)] $a(K[F])=-{\rm min}\left.\left\{\langle x_0,x\rangle\right
\vert\, x\in\mathbb{Z}^n\cap \{x\vert\, xB\leq b\}\right\}$.
\end{itemize}
\end{theorem}

\begin{proof} Let ${\mathcal H}=\{c_1,\ldots,c_r\}$.  
By duality \cite[Corollary~7.1a]{Schr}, we have the equality
\begin{equation}\label{aug31-06}
\mathbb{R}_+{\mathcal A}=H_{c_1}^+\cap\cdots\cap H_{c_r}^+.
\end{equation}
Observe that 
$\mathbb{R}_+{\mathcal A}\cap H_{c_i}$ is a proper face if $b_i=-1$ and it
is an improper face otherwise. 
From Eq.~(\ref{aug31-06}) we get that each facet of
$\mathbb{R}_+{\mathcal A}$ has the form  
$\mathbb{R}_+{\mathcal A}\cap H_{c_i}$ for some $i$. The relative
interior of the cone $\mathbb{R}_+{\mathcal A}$ is the union of its
facets.  Hence, using that $\mathcal{H}$ is an integral basis, 
we obtain the equality
\begin{equation}\label{subsa}
\mathbb{Z}^n\cap(\mathbb{R}_+{\mathcal A})^{\rm o}=\mathbb{Z}^n\cap
\{x\vert\,  
xB\leq b\}.
\end{equation}
Now, part (a) follows readily from Eqs.~(\ref{ejc8}) and
(\ref{subsa}). Part (b) follows from 
Eq.~(\ref{princeton-fall-07}) and part (a).  \end{proof}

Next we illustrate how to determine the canonical module and the
$a$-invariant using Theorem~\ref{canmod-tdi}.  

\begin{example} Let
$F=\{x_1,x_2,x_3,x_4,x_1x_2x_5,x_2x_3x_5,x_3x_4x_5,x_1x_4x_5\}$ and
let $\mathcal{A}$ be the set of exponent vectors of the monomials 
in $F$. Notice that $\mathcal{A}$ is a Hilbert basis and $\langle
x_0,v\rangle=1$ for 
$v\in \mathcal{A}$, where $x_0=(1,1,1,1,-1)$. 
An integral basis for $(\mathbb{R}_+\mathcal{A})^*$ is 
given by 
$$
\{e_1,e_2,e_3,e_4,e_5,(0,1,0,1,-1),(1,0,1,0,-1)\}.
$$
Then it is easy to verify that $\omega_{K[F]}$ is generated by the set
of all monomials $x^a$ such that $a=(a_i)$ is in the polyhedron $Q$
defined 
by the system:
$$
a_i\geq 1\, \forall\, i;\ \ a_1+a_3-a_5\geq 1;\ \ a_2+a_4-a_5\geq 1.
$$
The only vertex of the polyhedron $Q$ is $v_0=(1,1,1,1,1)$. Thus the
$a$-invariant of $K[F]$ is equal to $-\langle x_0,v_0\rangle=-3$.
\end{example}

\section{Canonical modules and integer rounding
properties}\label{section-on-canmod} 

In this section we give a description of the canonical module and the
$a$-invariant for subrings arising from systems with the integer
rounding property.

Let $A$ be a matrix of order $n\times q$ with entries in $\mathbb{N}$
such that $A$ has non-zero rows and non-zero columns. Let
$v_1,\ldots,v_q$ be the columns of $A$. For use 
below consider the set $w_1,\ldots,w_r$
of all $\alpha\in\mathbb{N}^n$ such that $\alpha\leq v_i$ for 
some $i$. Let $R=K[x_1,\ldots,x_n]$ be a polynomial ring over a field
$K$ and 
let 
$$
S=K[x^{w_1}t,\ldots,x^{w_r}t]\subset R[t]
$$
be the subring of $R[t]$ generated by
$x^{w_1}t,\ldots,x^{w_r}t$, where $t$ is a new variable. 
As $(w_i,1)$ lies in the hyperplane
$x_{n+1}=1$ 
for all $i$, $S$ is a standard $K$-algebra. Thus a monomial $x^at^b$
in $S$ has degree $b$. In what follows we assume that $S$ has this
grading. If $S$ is normal,
then according to Eq.~(\ref{ejc8}) the canonical module of $S$ is 
the ideal given by  
\begin{equation}\label{ejc9}
\omega_{S}=(\{x^at^b\vert\,
(a,b)\in \mathbb{N}{\mathcal B}\cap ({\mathbb R}_+{\mathcal B})^{\rm o}\}),
\end{equation}
where  ${\mathcal B}=\{(w_1,1),\ldots,(w_r,1)\}$ and 
$({\mathbb R}_+{\mathcal B})^{\rm o}$ is the interior of 
${\mathbb R}_+{\mathcal B}$ relative to 
${\rm aff}({\mathbb R}_+{\mathcal B})$, 
the affine hull of ${\mathbb R}_+{\mathcal B}$. In our case  
${\rm aff}({\mathbb R}_+{\mathcal B})=\mathbb{R}^{n+1}$. 

Let $\ell_0,\ell_1,\ldots,\ell_m$ be the vertices of 
$P=\{x\vert\, x\geq 0; xA\leq \mathbf{1}\}$, where $\ell_0=0$, and 
let $\ell_1,\ldots,\ell_p$ be the set of all maximal elements of 
$\ell_0,\ell_1,\ldots,\ell_m$ (maximal with respect to $\leq$). 

\begin{lemma}\label{may26-08} For each $1\leq i\leq p$ there is a
unique positive  
integer $d_i$ such that the non-zero entries of $(-d_i\ell_i,d_i)$ are
relatively prime. 
\end{lemma}

\begin{proof} If the non-zero rational entries of
$\ell_i$ are written in 
lowest terms, then $d_i$ is the least common multiple of the denominators. 
\end{proof}

\noindent {\it Notation} In what follows $\{\ell_1,\ldots,\ell_p\}$ is
the set of maximal elements of $\{\ell_0,\ldots,\ell_m\}$ 
and $d_1,\ldots,d_p$ are the
unique positive integers in Lemma~\ref{may26-08}.

\medskip

The next result complements a result of \cite{roundp}.

\begin{theorem}\label{can-mod-intr-norm} 
If the system $x\geq 0; xA\leq \mathbf{1}$ has the 
integer rounding property, then the subring 
$S=K[x^{w_1}t,\ldots,x^{w_r}t]$ is normal, the canonical module of
$S$ is given by   
\begin{equation}\label{march5-08-1}
\omega_S=\left(\left\{\left.x^at^b\right\vert\,
(a,b)\left(\begin{array}{rrrlrr} 
-d_1\ell_1&\cdots&-d_p\ell_p&e_1&\cdots& e_n\\
d_1        &\cdots& d_p     &\ 0  &  \cdots&0
\end{array}
\right)\geq\mathbf{1}\right\}\right),
\end{equation}
and the $a$-invariant of $S$ is equal to 
$-\max_i\{\lceil 1/d_i+|\ell_i|\rceil\}$. Here
$|\ell_i|=\langle\ell_i,\mathbf{1}\rangle$. 
\end{theorem}

\begin{proof} Note that in Eq.~(\ref{march5-08-1}) we regard 
$(-d_i\ell_i,d_i)$ and $e_j$ as column vectors for all $i,j$. 
The normality of $S$ follows 
from Theorem~\ref{round-up-char}. 
Recall that we have the following duality (see 
Section~\ref{i-r-p}):
\begin{eqnarray}
P=\{x\vert\, x\geq0;\langle x,w_i\rangle\leq{1}\, \forall i\}&=&{\rm 
conv}(\ell_0,\ell_1,\ldots,\ell_m),\nonumber\\
{\rm
conv}(w_1,\ldots,w_r)&=&
\{x\vert\, x\geq0;\langle x,\ell_i\rangle\leq{1}\forall
i\}=T(P),\label{equ-duality}
\end{eqnarray}
where $\{\ell_0,\ell_1,\ldots,\ell_m\}\subset\mathbb{Q}_+^n$ is the set of
vertices of $P$ and $\ell_0=0$. Therefore 
using Eq.~(\ref{equ-duality}) and the maximality of
$\ell_1,\ldots,\ell_p$ we obtain
\begin{equation}\label{equ-duality-1}
{\rm conv}(w_1,\ldots,w_r)=
\{x\vert\, x\geq0;\langle x,\ell_i\rangle\leq{1},\,\forall\,
i=1,\ldots,p\}. 
\end{equation}
We set $\mathcal{B}=\{(w_1,1),\ldots,(w_r,1)\}$. Note that 
$\mathbb{Z}\mathcal{B}=\mathbb{Z}^{n+1}$. 
From Eq.~(\ref{equ-duality-1}) it is seen that
\begin{equation}\label{march12-08}
\mathbb{R}_+\mathcal{B}=H_{e_1}^+\cap\cdots\cap H_{e_n}^+\cap
H_{(-d_1\ell_1,d_1)}^+\cap\cdots\cap H_{(-d_p\ell_p,d_p)}^+.
\end{equation}
Here $H_{a}^+$ denotes the closed halfspace 
$H_a^+=\{x\vert\, \langle
x,a\rangle\geq 0\}$ and $H_a$ stands for the hyperplane through the
origin with normal 
vector $a$. Notice that
$$
H_{e_1}\cap\mathbb{R}_+\mathcal{B},\ldots,H_{e_n}\cap
\mathbb{R}_+\mathcal{B},
H_{(-d_1\ell_1,d_1)}\cap\mathbb{R}_+\mathcal{B},\ldots,
H_{(-d_p\ell_p,d_p)}\cap\mathbb{R}_+\mathcal{B}
$$
are proper faces of
$\mathbb{R}_+\mathcal{B}$. Hence from Eq.~(\ref{march12-08}) 
we get that a vector $(a,b)$, with
$a\in\mathbb{Z}^n$, $b\in\mathbb{Z}$, is in the relative interior of 
$\mathbb{R}_+\mathcal{B}$ if and only if the entries of $a$ are
positive and $\langle(a,b),(-d_i\ell_i,d_i)\rangle\geq 1$ for all
$i$. Thus the 
required expression for $\omega_S$, i.e.,
Eq.~(\ref{march5-08-1}), follows using the normality of $S$
and the Danilov-Stanley formula given in Eq.~(\ref{ejc9}). 

It remains to prove the formula for $a(S)$, 
the $a$-invariant of $S$. Consider the vector 
$(\mathbf{1},b_0)$, where $b_0=\max_i\{\lceil 1/d_i+|\ell_i|\rceil\}$. 
Using Eq.~(\ref{march5-08-1}), it is not hard to see (by direct
substitution of $(\mathbf{1},b_0)$), that the monomial 
$x^{\mathbf{1}}t^{b_0}$ is in $\omega_S$. Thus from
Eq.~(\ref{princeton-fall-07}) we get $a(S)\geq -b_0$. Conversely if 
the monomial $x^at^b$ is in $\omega_S$, then again from
Eq.~(\ref{march5-08-1}) we get $\langle(-d_i\ell_i,d_i),(a,b)
\rangle\geq 1$ for all $i$ and $a_i\geq 1$ for all $i$, where $a=(a_i)$.
Hence 
$$
bd_i\geq 1+d_i\langle a,\ell_i\rangle\geq 1+d_i\langle
\mathbf{1},\ell_i\rangle=1+d_i|\ell_i|. 
$$
Since $b$ is an integer we obtain $b\geq \lceil 1/d_i+|\ell_i|\rceil$ for
all $i$. Therefore $b\geq b_0$, i.e., $\deg(x^at^b)=b\geq b_0$. As 
$x^at^b$ was an arbitrary monomial in $\omega_S$, by the formula 
for the $a$-invariant of $S$ given in Eq.~(\ref{princeton-fall-07}) 
we obtain that $a(S)\leq -b_0$. Altogether one has $a(S)=-b_0$, as
required.  \end{proof} 

A standard graded $K$-algebra $S$ is called {\it Gorenstein\/} if $S$ is
Cohen-Macaulay and $\omega_S$ is a principal ideal.

\begin{theorem}\label{may27-08-1} Assume that the system $x\geq
0${\rm ;} $xA\leq\mathbf{1}$ has 
the integer rounding property. If $S=K[x^{w_1}t,\ldots,x^{w_r}t]$ is
Gorenstein and $c_0=\max\{|\ell_i|\colon\, 1\leq i\leq p\}$ is an
integer, then $|\ell_k|=c_0$ for each $1\leq k\leq p$ such that
$\ell_k$ has integer entries.
\end{theorem}

\begin{proof} We proceed by contradiction. Assume that $|\ell_k|<c_0$ for
some integer $1\leq k\leq p$ such that $\ell_k$ is integral. We may assume that
$\ell_k=(1,\ldots,1,0,\ldots,0)$ and $|\ell_k|=s$. From
Eq.~(\ref{march12-08}) it follows that the monomial
$x^{\ell_k}t^{s-1}$ cannot be in $S$ because 
$(\ell_k,s-1)$ does not belong to $H_{(-d_k\ell_k,d_k)}^+$. Consider the 
monomial $x^at^b$, where $a=\ell_k+\mathbf{1}$, $b=b_0+s-1$ and 
$b_0=-a(S)$. We claim that the monomial $x^at^b$ is in $\omega_S$.  
By Theorem~\ref{can-mod-intr-norm} it suffices to show that 
$\langle(a,b),(-d_j\ell_j,d_j)\rangle\geq 1$ for $1\leq j\leq p$. 
Thus we need only show that $\langle(a,b),(-\ell_j,1)\rangle>0$ for
$1\leq j\leq p$. From the proof of Theorem~\ref{can-mod-intr-norm},
it is seen that $-a(S)=\max_i\{\lfloor|\ell_i|\rfloor\}+1$. Hence we
get $b_0=c_0+1$. One has the
following equalities
$$
\langle(a,b),(-\ell_j,1)\rangle=-|\ell_j|-\langle\ell_k,\ell_j\rangle+b_0+s-1=
-|\ell_j|-\langle\ell_k,\ell_j\rangle+c_0+s.
$$
Set
$\ell_j=(\ell_{j_1},\ldots,\ell_{jn})$. From Eq.~(\ref{march12-08}) we
get that the entries of each $\ell_j$ are less than or equal to
$1$. Case (I): If $\ell_{ji}<1$ for some $1\leq i\leq s$, then 
$s-\langle\ell_k,\ell_j\rangle>0$ and $c_0\geq |\ell_j|$. 
Case (II): $\ell_{ji}=1$ for $1\leq i\leq s$. Then $\ell_j\geq
\ell_k$. Thus by the maximality of $\ell_k$ we obtain
$\ell_j=\ell_k$. In both cases we obtain 
$\langle(a,b),(-\ell_j,1)\rangle>0$, as required. Hence the monomial 
$x^at^b$ is in $\omega_S$. Since $S$ is Gorenstein and $\omega_S$ is 
generated by $x^\mathbf{1}t^{b_0}$, we obtain that $x^at^b$ is a
multiple of $x^\mathbf{1}t^{b_0}$, i.e., $x^{\ell_k}t^{s-1}$ must be
in $S$, a contradiction. \end{proof}

\begin{theorem}\label{may27-08-2} Assume that the system $x\geq 0;
xA\leq \mathbf{1}$ 
has the 
integer rounding property. If $S=K[x^{w_1}t,\ldots,x^{w_r}t]$ 
and $-a(S)=1/d_i+|\ell_i|$ for $i=1,\ldots,p$, then $S$ is Gorenstein. 
\end{theorem}

\begin{proof} We set $b_0=-a(S)$ and
$\mathcal{B}=\{(w_1,1),\ldots,(w_r,1)\}$. The ring $S$ is normal by
Theorem~\ref{round-up-char}. Since the monomial
$x^{\mathbf{1}}t^{b_0}=x_1\cdots x_n t^{b_0}$ is in $\omega_S$, we
need only show that $\omega_S=(x^{\mathbf{1}}t^{b_0})$. Take
$x^at^b\in\omega_S$. It suffices to prove that 
$x^{a-\mathbf{1}}t^{b-b_0}$ is in $S$. 
Using Theorem~\ref{round-up-char}, one has 
$\mathbb{R}_+\mathcal{B}\cap\mathbb{Z}^{n+1}=
\mathbb{N}\mathcal{B}$. Thus we need only show that
the vector $(a-\mathbf{1},b-b_0)$ is in $\mathbb{R}_+\mathcal{B}$. From
Eq.~(\ref{march12-08}), the proof reduces to show 
that the vector $(a-\mathbf{1},b-b_0)$ is in $H_{(-\ell_i,1)}^+$ 
for $i=1,\ldots,p$.

As $(a,b)\in\omega_S$, from the description of $\omega_S$ given in 
Theorem~\ref{can-mod-intr-norm} we get
$$
\langle(a,b),(-d_i\ell_i,d_i) \rangle
=-\langle a,d_i\ell_i\rangle+bd_i\geq 1\ \Longrightarrow\ 
-\langle a,\ell_i\rangle\geq -b+1/d_i
$$
for $i=1,\ldots,p$. Therefore 
$$
\langle(a-\mathbf{1},b-b_0),(-\ell_i,1)\rangle=
-\langle a,\ell_i\rangle+|\ell_i|+b-b_0\geq
-b+{1}/{d_i}+|\ell_i|+b-b_0=0
$$
for all $i$, as required. \end{proof}

\begin{corollary}\label{may27-08-3} If $P=\{x\vert\, x\geq
0;xA\leq\mathbf{1}\}$ is 
an integral polytope, then the monomial subring
$S=K[x^{w_1}t,\ldots,x^{w_r}t]$  
is Gorenstein if and only if $a(S)=-(|\ell_i|+1)$ for $i=1,\ldots,p$. 
\end{corollary}

\begin{proof} Notice that if $P$ is integral, then $\ell_i$ has
entries in $\{0,1\}$ 
for $1\leq i\leq p$ and consequently $d_i=1$ for $1\leq i\leq p$. Thus
the result follows from Theorems~\ref{may27-08-1} and 
\ref{may27-08-2}. \end{proof}

\begin{example} Let $G$ be a pentagon with vertex set
$X=\{x_1,\ldots,x_5\}$, let $A$ be the incidence matrix of $G$ and let  
$$
S=K[t,x_1t,\ldots,x_5t,x_1x_2t,x_2x_3t,x_3x_4t,x_4x_5t,x_1x_5t].
$$
The system $x\geq 0;xA\leq \mathbf{1}$ has the integer rounding
property and the vertex set of $P=\{x\vert\, x \geq 0;xA\leq
\mathbf{1}\}$ is: 
$$
{\rm vert}(P)=\{0,\mathbf{1}/2,e_3+e_5,
e_2+e_5,e_2+e_4,e_1+e_4,e_1+e_3,e_1,e_2,e_3,e_4,e_5
\}.
$$
The maximal elements of ${\rm vert}(P)$ are
$$
\ell_1=\mathbf{1}/2,\ell_2=e_3+e_5,
\ell_3=e_2+e_5,\ell_4=e_2+e_4,\ell_5=e_1+e_4,\ell_6=e_1+e_3,
$$
$d_1=2$ and $d_i=1$ for $2\leq i\leq 6$. Notice that
$1/d_i+|\ell_i|=3$ for $i=1,\ldots,6$. Thus by
Theorems~\ref{can-mod-intr-norm} and \ref{may27-08-2}, the ring $S$ 
is Gorenstein and $a(S)=-3$. 
\end{example}

\begin{problem}\label{gorenstein-conjecture}\rm If $A$ is the
incidence 
matrix of a connected
graph and the system $x\geq 0$; $xA\leq\mathbf{1}$ has the integer
rounding
property, then the subring $S=K[x^{w_1}t,\ldots,x^{w_r}t]$ is
Gorenstein if and only if  
$-a(S)=1/d_i+|\ell_i|$ for $i=1,\ldots,p$.
\end{problem}

Note that the answer to this problem is positive if $A$ is the
incidence matrix of a bipartite graph because in this case $P$ is an
integral polytope and we may apply 
Corollary~\ref{may27-08-3}. If $A$ is the incidence matrix of a
connected non-bipartite graph $G$, E. Reyes has shown that $G$ is
unmixed if $S$ is Gorenstein. If $A$
is the incidence matrix of a graph, then it is seen that $d_i=1$ or
$d_i=1/2$ for each $i$. 

\subsubsection*{Subrings associated to the system $xA\leq \mathbf{1}$} 
Let $A$ be a matrix with entries in
$\mathbb{N}$ such that the system $xA\leq\mathbf{1}$ has integer
rounding property. As before we assume that the rows and
columns of $A$ are different from zero and that $v_1,\ldots,v_q$ are
the columns of $A$. In what follows we assume that $|v_i|=d$ 
for all $i$. 

The following lemma is not hard to show.

\begin{lemma}\label{jun8-08} If $|v_i|=d$ for all $i$. Then there are
isomorphisms
\begin{eqnarray*}
K[x^{v_1}t,\ldots,x^{v_q}t,t]\simeq
K[x^{v_1}t,\ldots,x^{v_q}t][T] \mbox{ and }\ \ \ \ \ 
\ \ \ \ \ \ \ \ \ \ \ \ & &  \\ 
K[x^{v_1}t,\ldots,x^{v_q}t]\simeq K[x^{v_1},\ldots,x^{v_q}]& &
\end{eqnarray*}
induced by $x^{v_i}t\mapsto x^{v_i}t$, $t\mapsto T$ and 
$x^{v_i}t\mapsto x^{v_i}$ respectively, where $T$ is a new
variable. 
\end{lemma}

Let $S$ be a homogeneous monomial subring and let $P_S$ be its toric
ideal. Recall that $S$ is called a {\it complete intersection\/} if
$P_S$ is a 
complete intersection, i.e., $P_S$ can be generated by ${\rm ht}(P_S)$
binomials, where ${\rm ht}(P_S)$ is the height of $P_S$. 
Let $c$ be a cycle of a graph $G$. A {\it chord\/} of $c$ is any edge
of $G$ joining two  
non adjacent vertices of $c$. A cycle without chords is called 
{\it primitive\/}.

\begin{proposition}\label{jun19-08} Let $G$ be a connected graph with $n$ 
vertices and $q$ edges and let $A$ be its incidence matrix. 
If the system $xA\leq\mathbf{1}$ has the integer rounding
property, then $K[x^{v_1}t,\ldots,x^{v_q}t,t]$ is 
a complete intersection if and only if $G$ is bipartite and 
the number of primitive cycles
of $G$ is equal to $q-n+1$. 
\end{proposition}

\begin{proof} $\Rightarrow$) By Corollary~\ref{jun8-08-1} the graph $G$ is
bipartite. From Lemma~\ref{jun8-08} it follows that 
$K[x^{v_1}t,\ldots,x^{v_q}t,t]$  is a complete intersection if and
only if  $K[G]=K[x^{v_1},\ldots,x^{v_q}]$ is a complete intersection.
Therefore by \cite{aron-jac} we get that $K[G]$ is a complete
intersection if and only if the number of primitive cycles of $G$ is
equal to $q-n+1$.

$\Leftarrow)$ By \cite{aron-jac} the ring $K[G]$ is a complete
intersection. Hence $K[x^{v_1}t,\ldots,x^{v_q}t,t]$ is a complete
intersection by Lemma~\ref{jun8-08}. \end{proof}

\noindent
{\bf Acknowledgments.} We gratefully acknowledge the computer algebra program
{\it Normaliz\/} \cite{B} which was invaluable in our
work on this paper. 
The third author also acknowledges the financial
support of  
CONACyT grant 49251-F and SNI. 

\bibliographystyle{plain}

\end{document}